\documentclass[12pt]{amsart}

\usepackage{latexsym}
\usepackage{amssymb}
\usepackage{amsmath}
\usepackage{amsfonts}

\newtheorem{theorem}{Theorem}
\newtheorem{lemma}{Lemma}
\newtheorem{corollary}{Corollary}

\newtheorem{theorem*}{Theorem}

\theoremstyle{definition}
\newtheorem{defn}[theorem]{Definition}

\newcommand{\bbZ}{\mathbb{Z}}
\newcommand{\bbN}{\mathbb{N}}

\newcommand{\ra}{\rangle}
\newcommand{\la}{\langle}

\begin{document}\title {On factoring of unlimited integers}
\author{Karel Hrbacek}
\address{Department of Mathematics\\                City College of CUNY\\                New York, NY 10031\\}                \email{khrbacek@icloud.com}
\date{April 30, 2023}

\begin{abstract}
Abdelmadjid Boudaoud asked whether every unlimited integer is a sum of a limited integer and a product of two unlimited integers. Assuming Dickson's Conjecture, the answer is negative.

The erroneous  proof of claim (2) on page 5 of the published version~\cite{H} is replaced by a correct one (see pages 6 - 7) and an answer to another question of Boudaoud and Bellaouar is given.
 \end{abstract}

\maketitle

Abdelmadjid Boudaoud \cite{B1, B2, B3} asked whether every large integer is close to a product of two large integers.
The question was made precise in the framework of nonstandard analysis by interpreting ``large'' as \emph{unlimited} (infinite, nonstandard) and ``close to'' as \emph{having a limited} (finite, standard) \emph{difference from}:

\bigskip
\emph{Is every unlimited integer a sum of a limited integer and a product of two unlimited integers? }

\bigskip
Symbolically: \emph{If $\omega \in \mathbb{Z}$ is unlimited, is 
\begin{equation} \tag{*}        \omega = s + \omega_1 \cdot \omega_2
\end{equation}
 where $s \in \mathbb{Z} $ is limited and $\omega_1, \omega_2 \in \mathbb{Z}$ are unlimited?}

\bigskip
We show that the question has negative answer assuming Dickson's Conjecture about primes in arithmetic progressions.
Following \cite{B1} we use the internal language of axiomatic nonstandard set theory such as IST or BST, but the argument works, with minor modifications, in any model-theoretic framework (ultraproducts, superstructures). Only the most basic ideas of nonstandard analysis are required, and those only for understanding of the conversion of the problem to an equivalent standard one.

\bigskip
First some simple observations.
Let $\omega$ be unlimited and $ \pi_1 \cdot \pi_2 \cdot \ldots \cdot \pi_\nu$ be the prime number decomposition of $|\omega|$. If at least two of the prime numbers are unlimited, or if $\nu$ is unlimited, then clearly $\omega = \omega_1 \cdot \omega_2$ for some unlimited integers $\omega_1, \omega_2$ and (*) holds with $s = 0$. 
Hence a counterexample to (*) has to have the form $\omega = a \cdot \pi$ where $\pi$ is an unlimited prime number and $a\in \mathbb{Z}$ is limited, $a \neq 0$. Without loss of generality we can assume that $\omega > 0$ and $\omega$ is a prime number. 
Indeed, if every unlimited prime number $\pi$ could be expressed in the form (*) as $\pi = s + \omega_1 \cdot \omega_2$, then $\omega = a \cdot \pi = (a\cdot s) + (a\cdot\omega_1) \cdot \omega_2$, which would have the required form (*).

\bigskip
If a prime number $\pi $ is a counterexample to (*),  then for each limited $s \in \mathbb{Z}$ there exist an unlimited prime number $\pi_s$ and a limited $a_s \in \mathbb{Z}$, $a_s > 0$, such that 
\begin{equation} \tag{**}  \pi - s = a_s \cdot \pi_s.
\end{equation}
Noticing that $a_0 = 1$ and $\pi = \pi_0$, we rewrite $(**)$ as 
\begin{equation} \tag{***}  a_s \cdot \pi_s -\pi_0 = - s.
\end{equation}
It follows that, for all limited positive integers $q, r$, the system of Diophantine equations 
$ a_s \cdot x_s - x_0 = -s, \quad 0 < |s| \le q$,
has a solution where all $x_s $ are prime numbers greater than $r$.
By Standardization, we can extend the external sequence $\langle a_s : s \in \mathbb{Z}, s \text{ limited} \rangle$ to a standard sequence $\langle a_s : s \in \mathbb{Z} \rangle$.
By Transfer we deduce that the existence of a counterexample to Boudaoud's question implies the following statement $\mathcal{S}$ of standard number theory:
 
\bigskip
\emph{There is a sequence $\langle a_s : s \in \mathbb{Z} \rangle$ such that $a_0 = 1$ and for all positive integers $q, r$ the system of Diophantine equations 
\begin{equation}\tag{$S_q$}     a_s \cdot x_s - x_0 = -s, \quad 0 < |s| \le q
\end{equation}
has a solution where all $x_s $, $|s| \le q$, are prime numbers greater than $r$.}

\bigskip
On the other hand, if there is a sequence as in the statement $\mathcal{S}$, then there is a standard one, by Transfer. Given such standard sequence $\langle a_s : s \in \mathbb{Z} \rangle$, we can take unlimited $q$ and $r$ and the corresponding solution $\langle x_s : |s| \le q \rangle$.  For each $s$  limited, $a_s$ is limited and $x_s$ is an unlimited prime number satisfying $x_0 - s = a_s \cdot x_s$; so $\omega = x_0$ is a counterexample to Boudaoud's question.
Therefore it suffices to construct a sequence $\langle a_s : s \in \mathbb{Z} \rangle$ as in $\mathcal{S}$.

\bigskip
The system  of equations $(S_q)$ has a solution  $\langle x_s : |s| \le q \rangle$ if and only if the system of congruences 
\begin{equation}\tag{$R_q$}  x \equiv s \mod  a_s, \quad 0 < |s| \le q
\end{equation}
has a solution $x_0$.
The obvious necessary condition for solvability of $(R_q)$ is
\begin{equation}\tag{$C$} \qquad (a_s,  a_t) \mid t - s 
\end{equation}
for all   $0< |s|, |t| \le q$.
It is an easy corollary to the Chinese Remainder Theorem that the condition $(C)$ is also sufficient for the existence of a solution to $(R_q)$; moreover, if $\bar{x}_0$ is one solution of  $(R_q)$, then every solution is given by $ x_0(k) = \bar{x}_0 + A^{q}  \cdot k$ where $A^{q}  = [a_{-q}, \ldots , a_q]$ is the least common multiple of $a_s,  |s| \le q$, and $k \in \mathbb{Z}$.
We let $A^{q}_s = A^q / a_s$; returning to the system $(S_q)$ we see that, assuming the condition $(C)$ is satisfied for all   $0< |s|, |t| \le q$, the system is solvable and all of its solutions are of the form  
\begin{equation}\tag{$F$} 
  x_s (k) = \bar{x}_s + A^{q}_s \cdot k, \quad |s| \le q   
\end{equation}
where $\langle \bar{x}_s :  |s| \le q \rangle$ is a particular solution of $(S_q)$ and $k \in \mathbb{Z}$.
We note that the solutions are given by a set of arithmetic progressions. 
To obtain the desired result, we need to show that there are solutions where all $x_s(k)$ are prime, for arbitrarily large $k$. 

\bigskip
\emph{Dickson's Conjecture} was formulated by Leonard Dickson in ~\cite{D}:

\bigskip
\emph{Let $\ell \ge 1$, $f_i(x) = a_i + b_i \cdot x$ with $a_i$ and $b_i$ integers, $b_i \ge 1$ (for $i = 1,\ldots, \ell$).
If there does not exist any integer $n > 1$  dividing all the products $\Pi_{i = 1}^{\ell} \; f_i(k)$, for every integer $k$, then there exist infinitely many natural numbers $m$ such that all numbers $f_1(m),\ldots, f_{\ell}(m)$ are prime.}

\bigskip
Dickson's Conjecture implies that in~$(F)$ there are arbitrarily large $k$ for which all $x_s(k)$, $|s| \le q$, are prime numbers,  provided the following congruence condition is satisfied:
$$(D) \quad \text{For every prime } p \text{ there is } k \text{ such that }  p \nmid x_{s}(k)  \text{ holds for all } |s| \le q.$$

So it suffices to show that for every prime $p$ there is a solution $\langle x_s : |s| \le   q \rangle$ of $(S_q)$ such that  $ p \nmid x_s $ holds for all $|s| \le q$.

\bigskip
The condition $(C)$, which guarantees solvability of $(R_q)$, and therefore of $(S_q)$, does not imply $(D)$ (consider the possibility $a_0 = a_1 =1, p=2$).
We formulate a condition that does, for all $q$.\\\\
$(E)$ \quad For every prime number $p$  and every $s \in \mathbb{Z}$:

If $a_s = p^n \cdot a'_s$ with $ p \nmid a'_s$,  then there is $r\in \mathbb{Z}$ such that  $p^{n+1} \mid a_r$  and  $r - s= u \cdot p^n$ 
 where $0 < u < p$.

\lemma
If the sequence $\langle a_s : s \in \mathbb{Z} \rangle$ satisfies the conditions (C) (for all $s,t \in \mathbb{Z}$) and (E), then for every $q > 0$ and every prime number $p$ the system $(S_q)$ has a solution $\langle x_s : |s| \le q \rangle$ such that  $p \nmid x_s $ holds for all $|s| \le q$.

\proof
We fix $q$ and $p$. Let $n$ be the highest exponent such that $p^{n} \mid a_s$ for some $|s| \le q$.
Let $q^* = q + p^{n+1}$. Let $\langle x_s : |s| \le q^* \rangle$ be a solution of the system $(S_{q^*})$.
Since the restriction of this solution to $|s| \le q$ is a solution of $(S_q)$, it suffices to prove that for this solution $p \nmid x_s $ holds for all $|s| \le q$.

We fix $s$ with $|s| \le q$, write  $a_s = p^n \cdot a'_s$ with $ p \nmid a'_s$, take $r$ as in $(E)$, and notice that $|r| \le q^*$. 
The equation $a_r \cdot x_r - a_s \cdot x_s = r -s$ follows from $(S_{q^*})$.
We thus have 
$p^{n+1} \cdot a'_r\cdot x_r - p^n \cdot a'_s \cdot x_s = u \cdot p^n$,
and after simplifying, 
$p \cdot a'_r\cdot x_r -  a'_s \cdot x_s = u$.
If $p \mid x_s$, then $p \mid u$, a contradiction with $0 < u < p$.
\qed

\bigskip
It remains to construct a sequence $\langle a_s : s \in \mathbb{Z} \rangle$ that satisfies $(C)$ and $(E)$.
We describe its terms $a_s$ by their prime factorization $ \Pi p^{n_p(s)}$.
The basic idea is to space those $a_s$ that are divisible by $p^n$ exactly $p^n$ steps apart. (This can be accomplished for $s \ge 0$ by simply putting $a_s = s+1$; however, we need a sequence that has this property and is defined  for all $s \in \mathbb{Z}$.)

\begin{defn}
For $p > 2$ we let the \emph{anchor} $$s(p,n) = (p^n + 1)/2;$$ we also let $$s(2,n) = (1 - (-2)^n )/3 = 
\Sigma_{i=0}^{n-1} (-2)^i.$$
\end{defn}

\lemma
If $m < n$, then $s(p,n ) \equiv s(p,m) \mod p^m$.

\proof
For $p > 2$ this is obvious from $s(p,n) - s(p,m) =  (p^n -p^m)/2 = p^m \cdot (p^{n-m} -1)/2$.
Also  $s(2,n) - s(2,m) = \Sigma_{i=m}^{n-1} (-2)^i
= \pm 2^m \cdot \Sigma_{i=0}^{{n-m -1}} (-2)^i $.
\qed

\bigskip
\begin{defn}
For  every  $s$, we let $n_p (s)$ be the highest exponent $n$ for which 
$s \equiv s(p,n) \mod p^n$.
A choice of anchors is \emph{admissible} if $n_p(s)$ exists, for all $p$ and~$s$.
\end{defn}
 
We note that, for $p > 2$,  $s(p,n) = (p^n + 1)/2 > 0$ and $s(p,n) - p^n= - (p^n -1)/2 < 0$. Hence for all $n$ such that  $|s| < (p^n -1)/2$ we have $s \not \equiv s(p,n) \mod p^n$, 
and $n_p(s)  \le \min \{ n : 2 |s| \le p^n \}.$
 Similarly, for all $n$ such that $|s| < |1 - (-2)^{n-1}|/3$ we have $s \not \equiv s(2,n) \mod 2^n$;
hence $n_2(s)  \le \min \{ n : 3 |s| \le 2^{n-1} \}.$
 These observations show that our choice of anchors is admissible. 
Taking $n = 1$ establishes that for $p > 2|s| + 1$ we have $n_p(s) = 0$, so the coefficients $a_s =  \Pi p^{n_p(s)}$ are well-defined. A table of $a_s$ for $|s| \le 12$ is computed below.

\begin{center}
\begin{tabular}{|c|c|c|c|c|c|c|c|c|c|c|c|c|}
\hline
$s$ & -12 &-11&-10&-9&-8&-7&-6&-5& -4 & -3 &-2 & -1\\
\hline
$a_s$ & $5^2$& $ 2 \cdot 23$ &  $3 \cdot 7$ & $ 2^2 \cdot 19$ &  17 & $2 \cdot 3 \cdot 5 $ & 13 & $2^4 \cdot 11$ & $3^3$ & $2 \cdot 7$& 5& $ 2^2 \cdot 3$\\
\hline
\end{tabular}
\end{center}

\begin{center}
\begin{tabular}{|c|c|c|c|c|c|c|c|c|c|c|c|c|c|}
\hline
$s$ &0 &1 &2 &3 &4 &5 &6 & 7 & 8 &9  & 10 &11 &12\\
\hline
$a_s$ &1 & 2 & 3 &$2^3 \cdot 5$ & 7 & $ 2 \cdot 3^2$& 11 &$2^2 \cdot 13$ &$ 3 \cdot 5$ &$2 \cdot 17$& 19& $ 2^5 \cdot 3\cdot 7$&23\\
\hline
\end{tabular}
\end{center}

\bigskip

\lemma
The sequence $\langle a_s : s \in \mathbb{Z} \rangle$ satisfies the condition~$(C)$  for all $s,t \in \mathbb{Z}$.
 
\proof
Fix $s, t \in \mathbb{Z}$ and a prime number $p$. By definition of $n_p(s)$,   $s \equiv s(p,n_p(s)) \mod p^{n_p(s)}$, and by definition of $a_s$,  $n_p (s)$ is the highest exponent $n$ for which $p^n \mid a_s$. Similarly, $t \equiv s(p,n_p(t)) \mod p^{n_p(t)}$, and  $n_p (t)$ is the highest exponent $n$ for which $p^n \mid a_t$. 
By Lemma 2, for $m = \min (n_p(s), n_p(t))$ we have 
$ s(p,n_p(s))  \equiv s(p,m)  \equiv s(p,n_p(t)) \mod p^{m}$; hence $s \equiv t \mod p^m$. 
Thus, for every prime $p$ the highest power of $p$ that divides both $a_s$ and $a_t$ also divides $t - s$.
From this the condition~$(C)$, to wit, $(a_s, a_t ) \mid t-s$, readily follows.

 \lemma
The sequence $\langle a_s : s \in \mathbb{Z} \rangle$ satisfies the condition~$(E)$.

\proof
Fix $s\in \mathbb{Z}$ and a prime number $p$ so that $a_s = p^n \cdot a'_s$ and $p \nmid a'_s$; by the construction of $a_s$, $n = n_p(s)$.
There is a unique $k \in \mathbb{Z}$ such that
$s(p, n+1)  +k \cdot p^{n+1} \le s < s(p, n+1)  + (k+1) \cdot p^{n+1} $; we let 
$r = s(p, n+1)  +k \cdot p^{n+1}$. Then $p^{n+1} \mid a_r$, so in particular $r \ne s$, and 
$ 0 < s- r < p^{n+1}$.
As  $s(p, n+1) \equiv s(p,n) \mod p^n$, we have $r \equiv s(p,n) \equiv s \mod p^n$.
Hence $s - r = u \cdot p^n$ and necessarily $0 < u < p$. This proves the condition~$(E)$.

\bigskip
We restate the final result in the language of model theory. 

\begin{theorem*}
Let $({}^*\mathbb{Z}, <, +, \times, 0, 1)$ be an elementary extension of $(\mathbb{Z}, <, +, \times, 0, 1)$ with $\mathbb{Z} \subset {}^*\mathbb{Z}$.
Assuming Dickson's Conjecture, there exist $\omega  \in {}^*\mathbb{Z} \setminus \mathbb{Z}$ such that
every integer in the galaxy of $\omega $, defined as $ \mathbf{G} (\omega) = \{ \omega - s \mid s \text{ limited}\}$, factors as $\omega - s = a_s \cdot \pi_s$ where $a_s \in \mathbb{Z}$ and $\pi_s$ is an unlimited prime number. 
In the example constructed above, the galaxy $\mathbf{G}(\omega)$ contains a unique prime number $\omega  = \pi_0$.
\end{theorem*}

\proof
The sequence $\langle a_s : s \in \mathbb{Z} \rangle$ constructed above is first-order definable in the language of  $(\mathbb{Z}, <, +, \times, 0, 1)$. The elementary extension assumption is sufficient to conclude that there are $q, r \in {}^*\mathbb{Z} \setminus \mathbb{Z}$ for which the system ($S_q$)
has a solution where all $x_s $ are prime numbers greater than $r$.

According to the construction, the coefficient $a_s$ is divisible by an odd prime $p $ if and only if  
$s = (p+1)/2 + k\cdot p$ for some $k \in   \mathbb{Z}$ if and only if
$2s - 1 = p\cdot(2k + 1)$ for some $k \in   \mathbb{Z}$.
Thus at least one such $p$ exists, except when $s = 1$ and $s = 0$.
The coefficient $a_1=2$.
We conclude that all  $\omega_s - s = a_s \cdot \pi_s \in \mathbf{G} (\omega)$ are composite except for $\omega_0 = \pi_0$.
\qed

\bigskip
The proofs that the sequence $\langle a_s : s \in \mathbb{Z} \rangle$ satisfies~$(C)$ and~$(E)$ go through for every admissible choice of anchors $s(p,n)$ such that, for every $s$, the exponents $n_p (s) = 0$ for all sufficiently large $p$. 
The anchors can be choosen so as to make the galaxy $\mathbf{G}(\omega)$ contain 

(1) no prime numbers or 

(2) infinitely many prime numbers.

 \bigskip
(1)  
The arguments leading to Theorem 1 are easily modified to allow relaxing the requirement $a_0 = 1$.
The following values of $s(p,n)$ guarantee that  $a_s \neq 1$ holds for all $s\in \mathbb{Z}$, and hence $\mathbf{G}(\omega)$ contains no prime numbers. 

Let $ p_1, p_2 , \ldots, p_k,\ldots $ be the increasing enumeration of odd primes.
We set $s(2,0) = 1$, $s(2,1) = 0$, $s(2,n) = 2^{n-1}$ for $n > 1$; $s(p_i, 0) = 1$, $s(p_{2i}, 1) = i$, 
$s(p_{2i+1}, 1) = -i$, and $s(p,n) = s(p,1) + p^{n-1}$ for all odd $p$ and all $n >1$.

It is easy to verify that this choice of anchors is admissible and that, for every $s$, the exponents $n_p(s) = 0$ for all sufficiently large $p$.

A table of $a_s$ for $|s| \le 12$ is partially  (up to $i = 5$, $p_{11} = 37$) computed below.
\begin{center}
\begin{tabular}{|c|c|c|c|c|c|c|c|c|}
\hline
$s$ & -12 &-11&-10&-9&-8&-7&-6&-5\\
\hline
$a_s$ & $2^3\cdot 3\cdot $ & $  \cdot  $ &  $2^2\cdot $ & $ 3 \cdot 5 \cdot 11 $ &  $2^4\cdot 7 $& $ \cdot  $ & $2^2\cdot3^2\cdot $  & $ 37  $\\
\hline
\end{tabular}
\end{center}

\begin{center}
\begin{tabular}{|c|c|c|c|c|c|c|c|c|c|}
\hline
  -4 & -3 &-2 & -1 &0 &1 &2 &3 &4 \\
\hline
 $2^3\cdot 5\cdot 29$ & $ 3\cdot 19$& $2^2\cdot 13$& $ 7 $&$2\cdot 3$ & 5 & $2^2\cdot 11$ &$3^2\cdot 17 $ & $2^3 \cdot 23$\\
\hline
\end{tabular}
\end{center}

\begin{center}
\begin{tabular}{|c|c|c|c|c|c|c|c|c|c|c|c|c|c|}
\hline
 5 &6 & 7 & 8 &9  & 10 &11 &12\\
\hline
$  31 $ & $ 2^2\cdot 3 \cdot 5 \cdot 7^2 $&$ \cdot  $ &$ 2^4 \cdot$ &$ 3^3\cdot 11^2 $& $2^2 \cdot$& $ 5 \cdot 13^2$& $2^3\cdot 3^2 \cdot $\\
\hline
\end{tabular}
\end{center}

\bigskip
(2) 
We construct values of $s(p,n)$ so that $\mathbf{G}(\omega)$ contains infinitely many prime numbers. 
(This replaces an erroneous argument in the published version.)

Let $p_0 =2, p_1=3, p_2,\ldots,p_k,\ldots$ enumerate prime numbers in the increasing order.
We construct a strictly increasing sequence $\la s_i \ra_{i \in \bbN}$ and the anchors $s(p_j,1)$ for $j \in \bbN$ by recursion so that
$s_i \not \equiv s(p_j,1) \mod p_j$ holds  for all $i,j \in \bbN$.
It follows that $a_{s_i} = 1 $ for all $i \in \bbN$.

Let $s_0 = 0$ and choose $s(p_0,1) \not \equiv s_0 \mod p_0$  ($s(p_0,1)  = 1$ will do).

At stage $k$ let $s_k$ be the least $s > s_{k-1}$ such that
$$s \not \equiv s(p_j,1) \mod p_j \;\text{  for }  0\le j \le k-1.$$

Then choose $s(p_k,1) < p_k$ so that 
$$s(p_k,1) \not \equiv s_i \mod p_i \;\text{  for }  0\le i \le k.$$

The choice of $s(p_k,1)$ is always possible because there are $k + 1$ values to be avoided,
 $p_k$ congruence classes $\mod p_k$  available,
and $k+1 < p_k$ for all $k\ge 0$.
Moreover, for $k \ge 4$ we have $p_k > 2k + 2$. 
It follows that the interval 
$(\frac{1}{4} p_k, \;  \frac{3}{4} p_k)$ contains more than $k+1$ values, so $s(p_k,1)$ can be chosen in this interval. 
If we also let  $s(2,n) = \Sigma_{i=0}^{n-1} (-2)^i$ and
$s(p_k, n+1) = s(p_k,1) + p_k^n$ for $k>0$, the choice of anchors is admissible and the exponents $n_p(s) = 0$ for all sufficiently large $p$.

The table showing the primes up to $p_7 = 19$ in the decomposition of each $a_s$ for $0 \le s \le 16$   is below. 
In the table we always choose the smallest value of $s(p_k, 1)$ that qualifies, bullets mark the values of $s_i$, and asterisks on a prime $p$ denote its appearance in the anchor ($s(p,1))$ position.

\begin{center}
\begin{tabular}{|c|c|c|c|c|c|c|c|c|}
\hline
$s$ & $0^\bullet$ &$1$&$2^\bullet$&3&4&5&$6^\bullet$&7\\
\hline
$a_s$ & $1 $ & $  2^*\cdot 3^*$ &  $1  $ & $  2^3\cdot 5^* \cdot 7^*    $ &  $ 3^2 $& $  2 $ & $ 1  $  & $2^2\cdot 3 \cdot 11^*  $\\
\hline
\end{tabular}
\end{center}

\begin{center}
\begin{tabular}{|c|c|c|c|c|c|c|c|c|c|}
\hline
  8 & 9 &10 &11 &$12^\bullet$ &13 &$14^\bullet$&15 &16 \\
\hline
 $5^2 $ & $2 \cdot 13^*\cdot 17^* \cdot 19^*$& $3^3\cdot 7^2 $& $ 2 ^4$&$1 $ &$2\cdot 3^2\cdot 5  $ & $ 1$ &$2 ^2 $ & $ 3 $\\
\hline
\end{tabular}
\end{center}

\bigskip\bigskip
\textbf{Addendum.}
Boudaoud and Bellaouar~\cite{B3} also asked:

\bigskip
\emph{If $\omega = s + \omega_1\cdot \omega_2$ with $s$ limited and $ \omega_1, \omega_2$ unlimited, is 
then also $\omega = s' + \omega'_1\cdot \omega'_2$ with $s'$ limited and $ \omega'_1 \sim \omega'_2$?}

\bigskip
We recall that $\omega_1 \sim \omega_2$ means that there is a limited $k \in \bbN$ such that $|\omega_1| \le k\cdot |\omega_2| $ and $|\omega_2| \le k\cdot |\omega_1| $, ie, the ratio 
$\omega_1/ \omega_2$ is appreciable.

\corollary
Assuming Dickson's conjecture, the answer is NO.

\proof
Let $({}^*\mathbb{Z}, <, +, \times, 0, 1)$ be an elementary extension of the structure $(\mathbb{Z}, <, +, \times, 0, 1)$ with $\mathbb{Z} \subset {}^*\mathbb{Z}$ and let $({}^{**}\mathbb{Z}, <, +, \times, 0, 1)$ be an elementary \emph{end extension} of $({}^*\mathbb{Z}, <, +, \times, 0, 1)$ with ${}^*\mathbb{Z} \subset {}^{**}\mathbb{Z}$. 

Fix a prime number $q\in {}^*\mathbb{Z}\setminus \bbZ$. The sequence  $\la a_s : s\in {}^*\mathbb{Z}\ra$
as in the proof of Theorem 1, with ${}^*\bbZ$ in place of $\bbZ$,  is definable in $({}^*\mathbb{Z}, <, +, \times, 0, 1)$. 
We modify it by choosing $s(q,1) = 1$ and $s(q, n) = 1 + q^{n-1}$ for $n > 1$. This makes  $q $ divide $a_1$ and hence $a_1\in {}^*\mathbb{Z}\setminus \bbZ$.
As in the proof of Theorem 1, there exist primes $\pi_s \in {}^{**}\mathbb{Z} \setminus {}^{*} \bbZ$ 
such that $\pi_0 - s = a_s \cdot \pi_s$ for all $s \in {}^*\bbZ$.
We let $\omega = \pi_0$; then $\omega  = 1 + a_1\cdot \pi_1$ where $ a_1, \pi_1$ are unlimited.
Assume that 
$\omega - s = a_s \cdot \pi_s= \omega_1 \cdot \omega_2$.
If  $\pi_s \; | \; \omega_1$, then $\omega_2\;|\; a_s$, 
hence $\omega_2 \in {}^* \bbZ$ and $\omega_1 \in {}^{**} \bbZ \setminus {}^* \bbZ$.
It follows that   $\omega_1 \not \sim \omega_2$. The case $\pi_s \; | \; \omega_2$ is similar.
\qed

\bigskip

\bibliographystyle{amsalpha}

\end{document}